\documentclass[12pt,a4paper]{article}
\usepackage{amssymb,amsmath,amsfonts,amscd}

\usepackage{url}

\reversemarginpar

\newcommand{\initpar}{\makebox[\parindent]{ }}
\newtheorem{twX}{Theorem}[section]
\newtheorem{lemX}{Lemma}
\newtheorem{propX}[twX]{Proposition}

\newtheorem{wnX}[twX]{Corollary}
\newtheorem{exaX}[twX]{Example}
\newtheorem{pytX}[twX]{Question}

\def\kropka{\hspace{-2mm}\mbox{\bf .}}

\def\kon{\mbox{$\square$}}

\newcommand{\btw}[1]{\begin{twX}\label{#1}\kropka}
\newcommand{\bbtw}[2]{\begin{twX}[#1]\label{#2}\kropka\hspace{2mm}}
\newcommand{\etw}{\end{twX}}

\newcommand{\blem}[1]{\begin{lemX}\label{#1}\kropka}
\newcommand{\bblem}[2]{\begin{lemX}[#1]\label{#2}\kropka\hspace{2mm}}
\newcommand{\elem}{\end{lemX}}

\newcommand{\bprop}[1]{\begin{propX}\label{#1}\kropka}
\newcommand{\bbprop}[2]{\begin{propX}[#1]\label{#2}\kropka\hspace{2mm}}
\newcommand{\eprop}{\end{propX}}

\newcommand{\bwn}[1]{\begin{wnX}\label{#1}\kropka}
\newcommand{\bbwn}[2]{\begin{wnX}[#1]\label{#2}\kropka\hspace{2mm}}
\newcommand{\ewn}{\end{wnX}}

\newcommand{\bexa}[1]{\begin{exaX}\label{#1}\kropka}
\newcommand{\bbexa}[2]{\begin{exaX}[#1]\label{#2}\kropka\hspace{2mm}}
\newcommand{\eexa}{\end{exaX}}

\newcommand{\bpyt}[1]{\begin{pytX}\label{#1}\kropka}
\newcommand{\bbpyt}[2]{\begin{pytX}[#1]\label{#2}\kropka\hspace{2mm}}
\newcommand{\epyt}{\end{pytX}}

\newenvironment{dowod}
 {\begin{bf}Proof. \end{bf}}{\kon \medskip \par}

\newcommand{\war}[1]{\initpar\mbox{$(#1)$}}
\newcommand{\impl}[2]{(#1) $\Rightarrow$ (#2)}

\newcommand{\f}{\varphi}

\begin{document}
\thispagestyle{empty}

\begin{center}
\begin{LARGE}

{\bf Strong divisibility and lcm-sequences}

\end{LARGE}

\bigskip

{\bf  Andrzej Nowicki}

\bigskip

Nicolaus Copernicus University,\\
 Faculty of Mathematics and
Computer Science,\\
 87-100 Toru\'n, Poland, (e-mail:  anow@mat.uni.torun.pl).


\bigskip


\end{center}

\medskip

\begin{abstract}
Let $R$ be a gcd-domain (for example let $R$ be a unique factorization domain),
and let $(a_n)_{n\geqslant1}$ be a sequence of nonzero elements in $R$.
We prove that
$\gcd(a_n,a_m)=a_{\gcd(n,m)}$ for all $n,m\geqslant1$
if and only if
$$a_n=\prod\limits_{d\mid n} c_d\quad\mbox{for} \ n\geqslant1,
$$
where
$c_1=a_1$ and $c_n=\mbox{lcm}(a_1,a_2,\dots,a_n)/\mbox{lcm}(a_1,a_2,\dots,a_{n-1})$
for $n\geqslant2$.
All equalities with gcd and lcm are determined up to units of $R$.
\end{abstract}






\section*{Introduction}
\initpar
Let  $R$ be a commutative unique factorization domain,
and let $(a_n)_{n\geqslant1}$ be a sequence of nonzero elements of $R$.
Put $c_1=a_1$  and
$c_n=\mbox{lcm}(a_1,a_2,\dots, a_n)/\mbox{lcm}(a_1,a_2,\dots, a_{n-1})$
for  $n\geqslant2$.
Then $(c_n)$ is a sequence of nonzero elements of $R$ and it is well-defined
up to units.
We will say that $(c_n)$ is the {\it lcm-sequence of $(a_n)$}.
A sequence  $(a_n)_{n\geqslant1}$
of nonzero elements of $R$ is called a {\it strong divisibility sequence}
if $a_{\gcd(m,n)}$ is a greatest common divisor of $a_m$ and $a_n$ for all $m,n\geqslant1$.

\smallskip

Consider the case  $R=\Bbb Z[x]$, and let $a_n=x^n-1$ for $n\geqslant1$.
The $n$th cyclotomic polynomial is commonly defined by the formula
$\Phi_n(x)=(x-\omega_1)\cdots(x-\omega_{\varphi(n)})$,
where $\f$ is the Euler totient function and $\omega_1,\dots,\omega_{\f(n)}$
are the primitive $n$th roots of unity.
It is well known that for any $n$ we have the equality
$
a_n=\prod_{d\mid n} \Phi_d(x).
$
Applying the M\"obius inversion formula, we obtain the equality
$
\Phi_n(x)=\prod_{d\mid n} a_n^{\mu(n/d)},
$
in which  $\mu$ is the M\"obius function.
This formula may be used as a definition of $\Phi_n(x)$.
Recently, Tomasz Ordowski presented a new alternative
definition of $\Phi_n(x)$. He
observed that $\left(\Phi_n(x)\right)$ is the lcm-sequence of $(a_n)$,
that is, $\Phi_1(x)=x-1$ and
$$
\Phi_n(x)=\frac{\mbox{lcm}\Big(x^1-1, x^2-1,\dots,x^n-1\Big)}
{\mbox{lcm}\Big(x^1-1,x^2-1,\dots,x^{n-1}-1\Big)}
\quad\mbox{for}\quad n\geqslant2.
$$
Tomasz Ordowski observed also that there exists a similar formula
for integer values of $\Phi_n(x)$.
If $b\geqslant2$ is an integer, then $\Phi_1(b)=b-1$, and
$$
\Phi_n(b)=\frac{\mbox{lcm}\Big(b^1-1,b^2-1,\dots,b^n-1\Big)}
{\mbox{lcm}\Big(b^1-1,b^2-1,\dots,b^{n-1}-1\Big)}
\quad\mbox{for}\quad n\geqslant2.
$$
This means that the sequence $\left(\Phi_n(b)\right)_{n\geqslant1}$
is equal to the lcm-sequence of $\left(b^n-1\right)_{n\geqslant1}$.

\medskip

It is well known that $\left(x^n-1\right)_{n\geqslant1}$
and $\left(b^n-1\right)_{n\geqslant1}$ are strong divisibility sequences
in $\Bbb Z[x]$ and $\Bbb Z$, respectively.
In this article we show that
every strong divisibility sequence
has the above properties.

\medskip

In this article a sequence  $(a_n)_{n\geqslant1}$
of nonzero elements of $R$ is called  {\it special}
if there exists a sequence $(b_n)_{n\geqslant1}$ of
nonzero elements in $R$ such that
$a_n=\prod_{d\mid n} b_d$ for any $n\geqslant 1$.
Note that if a sequence $(a_n)$ is special, then
its associated sequence $(b_n)$ is unique.
The uniqueness  follows
from the M\"obius inversion formula,
applied to the field of fractions of $R$; we have
$b_n=\prod_{d\mid n} a_d^{\mu(n/d)}$ for any $n\geqslant1$.

It is known (Mathematical Olympiad, Iran 2001, Problem 5,
see for example \cite{AAA} p. 218) that
every  strong divisibility sequence in $\Bbb Z$ is special.
Recently, Nathan Bliss, Ben Fulan, Stephen Lovett and Jeff Sommars
proved in \cite{Bei} that the same is true for strong divisibility sequences
in arbitrary unique factorization domains.
Thus, if $(a_n)$ is a strong divisibility sequence in a  factorization domain $R$,
then there exists a unique sequence $(b_n)$ of nonzero elements in $R$
with $b_1=a_1$ and $a_n=\prod_{d\mid n} b_d$.
The authors of \cite{Bei} proved also that there exist
special sequences which are not strong divisibility sequences.

\medskip

The main theorem of this article states that
if $(a_n)_{n\geqslant1}$ is a sequence of nonzero elements in $R$
and $(c_n)$ is the lcm-sequence of $(a_n)$,
then
$(a_n)$ is a strong
divisibility sequence if and only if
$$
a_n=\prod_{d\mid n} c_d\quad\mbox{for all}\quad n\geqslant1.
$$
Applying the main theorem to the sequences $(x^n-1)$ and $(b^n-1)$,
we obtain the above mentioned formulas of Ordowski.
It is well-known that the Fibonacci sequence, the sequences of
the form $(a^n-b^n)$ with $\gcd(a,b)=1$,
and other known sequences have the strong divisibility property.
Applying the main theorem and the  M\"obius inversion formula
to these sequences, we obtain new formulas for some lcm-sequences
and strong divisibility sequences.

\medskip

We restricted our attention to an arbitrary unique factorization domain $R$.
The main theorem of this article we prove for  more general domains.
We assume only that $R$ is a gcd-domain.

\section{Notations and preliminary facts}

Let $S$ be a domain, that is, $S$ is a commutative ring with identity
without zero divisors. Let $a_1,\dots,a_n$ be nonzero elements  of $S$.
A {\it greatest common divisors} (abbreviated as gcd)
of $a_1,\dots,a_n$
is a nonzero element $d$ in $S$ such that
$d\mid a_i$ for $i=1,\dots,n$, and
if $0\neq d'\in S$, $d'\mid a_i$ for $i=1,\dots,n$, then $d'\mid d$.
A {\it least common multiple} (abbreviated as lcm)
of $a_1,\dots,a_n$
is a nonzero element $m$ in $S$ such that
$a_i\mid m$ for $i=1,\dots,n$, and
if $0\neq m'\in S$, $a_i \mid m'$ for $i=1,\dots,n$, then $m\mid m'$.

\medskip

Throughout this article we assume that $R$ is a gcd-domain,
that is, $R$ is a domain and any two nonzero elements in $R$
have a greatest common divisor (\cite{Kapl}).
It follows from the definition of gcd-domains, that if $n\geqslant1$ and $a_1,\dots,a_n$ are nonzero
elements in $R$, then there exist a greatest common divisor and a least
common multiple of $a_1,\dots,a_n$. We adopt the notation
$(a_1,\dots,a_n)$ and $[a_1,\dots,a_n]$ for the
greatest common divisor and the least common multiple,
respectively, of $a_1,\dots,a_n$. Note that the elements $(a_1,\dots,a_n)$
and $[a_1,\dots,a_n]$ are determined only up to units.
We can allow this ambiguity in our article.

Unique factorization domains,
Bezout domains and valuation domains belong to the class of gcd-domains. 
Every gcd-domain is integrally closed. If $R$ is a gcd-domain, then the polynomial ring
$R[x]$ is also a gcd-domain (\cite{Kapl}).

We say that two nonzero elements  $a,b\in R$  are {\it relatively prime},
if $(a,b)=1$.   It should be carefully noted
that we are not assuming that $(a,b)$
is a linear combination of $a$ and~$b$.
Let us note some properties of gcd and lcm.

\bprop{Prop}
If $a,b,c$ are nonzero elements of a gcd-domain $R$, then:\smallskip\\
\war{1}  $((a,b),c)=(a,(b,c))=(a,b,c)$, \  $[[a,b],c]=[a,[b,c]]=[a,b,c]$;\smallskip\\
\war{2} $(ac,bc)=(a,b)c$, \quad $[ac,bc]=[a,b]c$;\smallskip\\
\war{3} $(a,b)[a,b]=ab$;\smallskip\\
\war{4} if $d=(a,b)$, then $(a/d,b/d)=1$;\smallskip\\
\war{5} if $(a,b)=1$ and $(a,c)=1$, then $(a,bc)=1$;\smallskip\\
\war{6} if $(a,bc)=1$, then  $(a,b)=1$ and  $(a,c)=1$;\smallskip\\
\war{7} if $(a,b)=1$ and $a\mid bc$, then  $a\mid c$;\smallskip\\
\war{8} if $(a,b)=1$, $a\mid c$ and $b\mid c$, then  $ab\mid c$.
\eprop

Recall that the equalities which appear in this proposition
are determined up to units.
The proof of this proposition is elementary, so we omit it.

\medskip

Let $(a_n)_{n\geqslant1}$ be a sequence of nonzero elements of
a gcd-domain $R$.
 Recall that $(a_n)$  is called a {\it strong divisibility sequence}
if $\left(a_m,a_n\right)=a_{(m,n)}$
for all  $m,n\geqslant1$.
The following two theorem from \cite{Bei} will play an important role.
They are stated in \cite{Bei} for unique factorization domains, but
they are also valid (with the same proofs) for arbitrary gcd-domains.

\bbtw{\cite{Bei}}{TwA}
If $(a_n)_{n\geqslant1}$ is a strong divisibility sequence of a gcd-domain
$R$, then there exists a unique sequence $(b_n)_{n\geqslant1}$
of nonzero elements in $R$ such that $a_n=\prod_{d\mid n}b_d$
for all $n\geqslant1$.
\etw

\bbtw{\cite{Bei}}{TwB}
Let  $R$ be a  gcd-domain
and let $(b_n)_{n\geqslant1}$ be a sequence of nonzero elements of $R$.
Let $(a_n)_{n\geqslant1}$ be the sequence such that $a_n=\prod_{d\mid n} b_d$.
The sequence $(a_n)$ is a strong divisibility sequence if and only if,
for all positive integers $m$ and $n$ such that $m\nmid n$ and $n\nmid m$,
the elements $b_m$ and $b_n$ are relatively prime.
\etw

Assume again that $(a_n)_{n\geqslant1}$ is a sequence
of nonzero elements of a gcd-domain $R$.
Put
$$e_1=1\quad\mbox{and}\quad  e_{n+1}=[e_n,a_n]\quad\mbox{for } \ n\geqslant 1.
$$
Then $e_{n+1}=[a_1,\dots,a_n]$ and $e_n\mid e_{n+1}$
for all $n\geqslant 1$.
Denote by $c_n$ the element $e_{n+1}/e_n$. Then  $(c_n)_{n\geqslant1}$  is a sequence
of nonzero elements in $R$. We call it the
{\it lcm-sequence} of $(a_n)$. Note, that $c_1=a_1$ and
$$
c_n=\frac{[a_1,\dots,a_n]}{[a_1,\dots,a_{n-1}]}\quad
\mbox{for all } \ n\geqslant2.
$$

If for example $(a_n)$ is a constant sequence,
$a_n=a$ for all $n\geqslant 1$ with $0\neq a\in R$, then
$c_1=a$ and $c_n=1$ for all $n\geqslant2$.
The lcm-sequence of the geometric sequence $a_n=q^n$ with $0\neq q\in R$
is equal to the constant sequence $c_n=q$.
If $R=\Bbb Z$ and $a_n=n!$, then $c_n=n$.
The lcm-sequence of the sequence $a_n=n$ is of the form $(c_n)$, where
$$
c_n=\left\{
\begin{array}{ll}
p,&\mbox{if } \ n=p^s \ \mbox{ for some prime } p \mbox{ and integer } \ s\geqslant1,\smallskip\\
1,&\mbox{otherwise}.
\end{array}\right.
$$
The first few terms of the lcm-sequence of  the
triangular numbers
$t_n=n(n+1)/2$ are:
$
1,\ 3,\ 2,\ 5,\ 1,\ 7,\ 2,\ 3,\ 1,\ 11,\ 1,\ 13,\ 1,\ 1,\ 2,\ 17,\ 1,\ 19.
$

\section{The main theorem}

\btw{TwG}
Let $(a_n)_{n\geqslant1}$ is a sequence of nonzero elements
in a gcd-domain $R$, and let $(c_n)_{n\geqslant1}$ be the lcm-sequence
of $(a_n)$.
Then the following two conditions are equivalent.\smallskip\\
\war{1}
$(a_n)$ is a strong divisibility sequence.\smallskip\\
\war{2}
$a_n=\prod_{d\mid n} c_d$ for all $n\geqslant1$.
\etw

\begin{dowod}
\impl{1}{2}.
Assume that $(a_n)$ is a strong divisibility sequence.
Then, by Theorem \ref{TwA}, there exists a unique sequence $(b_n)_{n\geqslant1}$
of nonzero elements in $R$, such that $a_n=\prod_{d\mid n} b_d$
for all $n\geqslant1$.
Consider the sequence $(e_n)_{n\geqslant1}$ defined by  $e_1=1$ and
$e_{n+1}=[e_n,a_n]$ for $n\geqslant1$.
We will prove, using an induction with respect to $n$,
that for every $n\geqslant1$ we have the equality
$$
e_{n+1}=b_1b_2\cdots b_n.
\leqno{(\ast)}
$$
For $n=1$ it is obvious, because $e_2=[e_1,a_1]=a_1=b_1$.
Now let $n\geqslant2$, and assume that $e_{n}=b_1\cdots b_{n-1}$.
Then we have
$$
e_{n+1}=[e_{n},a_{n}]
=\left[\prod\limits_{k=1}^{n-1}b_k, \
\prod\limits_{d\mid n}b_d\right]
=[A\cdot B, \ A \cdot b_n],
$$
where
$A$ is the product of the elements $b_d$ with  $d<n$ and $d\mid n$,
and  $B$ is the product of all elements of the
form $b_d$
with  $d<n$ and $d\nmid n$.
Since $(a_n)$ is a strong divisibility sequence, we know
by Theorem \ref{TwB} and Proposition \ref{Prop} that the elements
$B$ and $b_n$ are relatively prime.
Hence, $[B, \ b_n]=B b_n$ and so, we have:
$$
e_{n+1}=[A B, \ A b_n]
=A[B, \ b_n]
=A B b_n=\left(\prod_{k=1}^{n-1}b_k\right)\cdot b_n
=b_1\cdots b_n.
$$
Thus, by induction, we have  a proof of $(\ast)$.
But, by the definition of the lcm-sequence of $(a_n)$,  $c_1=a_1$
and for $n\geqslant 2$ we have
$c_n=e_{n+1}/e_n=(b_1\cdots b_n)/(b_1\cdots b_{n-1})=b_n$.
Whence $b_n=c_n$ for all $n\geqslant1$.
Therefore, $a_n=\prod_{d\mid n} b_d=\prod_{d\mid n} c_d$ for $n\geqslant1$.

\medskip

\impl{2}{1}.
Assume now that $a_n=\prod_{d\mid n} c_d$ for $n\geqslant1$,
and let $m\geqslant 1$ be a fixed integer.
Denote by $U$
the product of the elements $c_d$ with  $d<m$ and $d\mid m$,
and  denote by $V$  the product of all elements of the
form $c_d$
with  $d<m$ and $d\nmid m$.
Then we have:
$$
UVc_m=\prod_{i=1}^m c_i=e_{m+1}=[e_m,a_m]
=\left[\prod_{i=1}^{m-1} c_i, \ \prod_{d\mid m} c_d\right]
=\left[UV, \ Uc_m\right]=U[V,c_m].
$$
Hence,  $[V,c_m]=Vc_m$ and hence,
the elements $V$ and $c_m$ are relatively prime.
Recall that  $V$  is the product of all elements
$c_d$ with  $d<m$ and $d\nmid m$.
This implies (by Proposition~\ref{Prop}) that if $d<m$ and $d\nmid m$,
then the elements $c_d$ and $c_m$ are relatively prime.

We proved that if $n,m\geqslant1$, $n\mid m$ and $m\nmid n$,
then the elements $c_n$ and $c_m$ are relatively prime.
Therefore, by Theorem \ref{TwB}, the sequence $(a_n)$ is a strong divisibility sequence.
\end{dowod}

It is well-known that the sequence $a_n=x^n-1$ is a strong
divisibility sequence in the polynomial ring $\Bbb Z[x]$. In this case
we have the equalities $x^n-1=\prod_{d\mid n}\Phi_d(x)$, where
each $\Phi_n(x)$ is the $n$th cyclotomic polynomial.
We know also, by Theorem \ref{TwG}, that $x^n-1=\prod_{d\mid n} c_d$,
where $(c_n)_{n\geqslant1}$ is the lcm-sequence of $(a_n)$.
Applying the M\"obius inversion formula we obtain that 
$\Phi_n(x)=c_n$ for all $n\geqslant1$, and hence, we have: 

\bbwn{Ordowski}{WnA} Let $a_n=x^n-1$. Then $\Phi_1(x)=a_1$ and
$$\displaystyle
\Phi_n(x)=\frac{\Big[a_1,a_2,\dots,a_{n-1}, a_n\Big]}
{\Big[a_1,a_2,\dots,a_{n-1}\Big]}
\quad\mbox{for}\quad n\geqslant2.
$$
\ewn

By the same way we obtain similar formulas for integer
values of cyclotomic polynomials. Let $b\geqslant2$ be an integer
and let $a_n=b^n-1$ for $n\geqslant1$. Then $(a_n)_{n\geqslant1}$
is a strong divisible sequence in $\Bbb Z$ and
$b^n-1=\prod_{d\mid n} \Phi_d(b)$ for $n\geqslant1$.
Applying Theorem \ref{TwG} and the M\"obius inversion formula
we again obtain that
$\left(\Phi_n(b)\right)_{n\geqslant1}$ is the lcm-sequence of $(a_n)$.
Therefore we have the following corollary.

\bbwn{Ordowski}{WnB} Let $b\geqslant2$ be an integer and let $a_n=b^n-1$ for $n\geqslant 1$.
Then
$$\displaystyle
\Phi_n(b)=\frac{\Big[a_1,a_2,\dots,a_{n-1}, a_n\Big]}
{\Big[a_1,a_2,\dots,a_{n-1}\Big]}
\quad\mbox{for}\quad n\geqslant2.
$$
\ewn

In particular, $\left[M_1,\dots,M_n\right]
/\left[M_1,\dots,M_{n-1}\right]=\Phi_n(2)$, where each $M_n$ is the
Mersenne number $2^n-1$.

It is easy to prove and it is well-known
that the sequence $a_n=x^n-y^n$ is a strong divisibility
sequence in the polynomial ring $\Bbb Z[x,y]$. In this case we have the
equalities $x^n-y^n=\prod_{d\mid n} \Psi_n(x,y)$,
where $\Psi_n(x,y)$ are polynomials in $\Bbb Z[x,y]$ defined by
$$
\Psi_n(x,y)=y^{\f(n)}\Phi_n(x/y)=\prod_{d\mid n}\left(x^d-y^d\right)^{\mu(n/d)}
$$
for all $n\geqslant 1$ (see for example \cite{Mot7}).
Thus, $[a_1,\dots,a_n]/[a_1,\dots,a_{n-1}]=\Psi_n(x,y)$
for all $n\geqslant2$.
Similarly, if $u>v\geqslant1 $ are relatively prime integers,
then the sequence $a_n=u^n-v^n$ is a strong divisible sequence in
$\Bbb Z$,
and again we have the equalities
$[a_1,\dots,a_n]/[a_1,\dots,a_{n-1}]=\Psi_n(u,v)$
for all $n\geqslant2$.

Examples and properties of strong divisibility sequences
can be found in many books and
articles (see for example:  \cite{Kim},  \cite{Sch}, \cite{Bei}, \cite{MaN}, \cite{Now} 150-152).
For each such sequence we have, by Theorem \ref{TwG},
an associated  formula connected with its lcm-sequence.
Consider several examples.
The sequence of repunits,
$$
u_n=\underbrace{11\dots1}_n=\frac{10^n-1}{9},
$$
is a strong divisible sequence in $\Bbb Z$.
The first few terms of the lcm-sequence are:
$$
1,\ u_2, \ u_3, \ 101, u_5, \ 91, \ u_7, \  10001, \   1001001, \  9091, \  u_{11}, \  9901, \ u_{13}, \  909091, \  90090991.
$$

The sequence of the Fibonacci numbers,  $F_1=F_2=1,  \ F_{n+2}=F_{n+1}+F_n$, is a strong divisibility sequence in $\Bbb Z$ (\cite{Vor}). 
The first few terms of its lcm-sequence are:
$$
1,\ 1,\ 2,\ 3,\ 5,\ 4,\ 13,\ 7,\ 17, \  11,\  89,\  6,\  233,\  29,\  61,\  47,\  1597,\  19,\  4181,\  41,\  421,\  199.
$$
Let 
$\left(F_n(x)\right)_{n\geqslant1}$ be the sequence of Fibonacci polynomials, that is, 
$$
F_1(x)=1, \ F_2(x)=x\quad\mbox{and } \ F_{n+2}(x)=xF_{n+1}(x)+F_n(x)\quad\mbox{for } \  n\geqslant1. 
$$
It is a strong divisibility sequence in $\Bbb Z[x]$ (\cite{Hog}).
The Chebyshev polynomials of the second kind (\cite{Bar}) ,
which satisfy $U_0(x)=1$,  $U_1(x)=2x$,  $U_{n+2}(x)=2xU_{n+1}-U_n(x)$, also form a strong divisibility sequence in $\Bbb Z[x]$ (\cite{Hog}). 
The same holds for the $3$-variable polynomials $S_n(x,y,z)$  which satisfy $S_1=1$, $S_2=x$,  $S_n=xS_{n-1}+yS_{n-2}$ for $n$ even and $S_n=zS_{n-1}+yS_{n-2}$ for $n$ odd (\cite{Hog}. 
Also, other examples of strong divisibility polynomial sequences are constructed and discussed in \cite{Bei}.

\medskip

Let us end this article with the following corollary, which
is an immediately consequence of Theorem \ref{TwG} and
Theorem 3 in \cite{Bei}.

\bwn{WnBei}
Let $(a_n)_{n\geqslant1}$ be a strong divisibility sequence
in a gcd-domain $R$.
If $n\geqslant2$ has prime factorization $n=p_1^{\alpha_1}\cdots p_s^{\alpha_s}$,
then
$$
\frac{\left[a_1,\dots,a_n\right]}{\left[a_1,\dots,a_{n-1}\right]}
=\frac{a_n}{\left[a_{n/p_1},a_{n/p_1},\dots,a_{n/p_s}\right]}\ .
$$
\ewn

\medskip

{\bf Acknowledgment.}
The author thanks  Mr. Tomasz Ordowski for his interesting formulas.


\end{document}